\documentclass[sw20lart]{article}
\usepackage{amsfonts}

\makeatletter
\newcount\@hour\newcount\@minute\chardef\@x10\chardef\@xv60
\def\tcitime{
\def\@time{%
  \@minute\time\@hour\@minute\divide\@hour\@xv
  \ifnum\@hour<\@x 0\fi\the\@hour:%
  \multiply\@hour\@xv\advance\@minute-\@hour
  \ifnum\@minute<\@x 0\fi\the\@minute
  }}%

\@ifundefined{hyperref}{}{}

\@ifundefined{qExtProgCall}{\def\qExtProgCall#1#2#3#4#5#6{\relax}}{}
\def\QCTOpt[#1]#2{%
  \def\QCTOptB{#1}
  \def\QCTOptA{#2}
}
\def\QCTNOpt#1{%
  \def\QCTOptA{#1}
  \let\QCTOptB\empty
}
\def\Qct{%
  \@ifnextchar[{%
    \QCTOpt}{\QCTNOpt}
}
\def\QCBOpt[#1]#2{%
  \def\QCBOptB{#1}
  \def\QCBOptA{#2}
}
\def\QCBNOpt#1{%
  \def\QCBOptA{#1}
  \let\QCBOptB\empty
}
\def\Qcb{%
  \@ifnextchar[{%
    \QCBOpt}{\QCBNOpt}
}
\def\PrepCapArgs{%
  \ifx\QCBOptA\empty
    \ifx\QCTOptA\empty
      {}%
    \else
      \ifx\QCTOptB\empty
        {\QCTOptA}%
      \else
        [\QCTOptB]{\QCTOptA}%
      \fi
    \fi
  \else
    \ifx\QCBOptA\empty
      {}%
    \else
      \ifx\QCBOptB\empty
        {\QCBOptA}%
      \else
        [\QCBOptB]{\QCBOptA}%
      \fi
    \fi
  \fi
}
\newcount\GRAPHICSTYPE
\GRAPHICSTYPE=\z@
\def\GRAPHICSPS#1{%
 \ifcase\GRAPHICSTYPE
   \special{ps: #1}%
 \or
   \special{language "PS", include "#1"}%
 \fi
}%
\def\graffile#1#2#3#4{%
    \leavevmode
    \raise -#4 \BOXTHEFRAME{%
        \hbox to #2{\raise #3\hbox to #2{\null #1\hfil}}}%
}%
\def\draftbox#1#2#3#4{%
 \leavevmode\raise -#4 \hbox{%
  \frame{\rlap{\protect\tiny #1}\hbox to #2%
   {\vrule height#3 width\z@ depth\z@\hfil}%
  }%
 }%
}%
\newcount\draft
\draft=\z@

\newif\ifwasdraft
\wasdraftfalse

\def\GRAPHIC#1#2#3#4#5{%
 \ifnum\draft=\@ne\draftbox{#2}{#3}{#4}{#5}%
  \else\graffile{#1}{#3}{#4}{#5}%
  \fi
 }%
\def\addtoLaTeXparams#1{%
    \edef\LaTeXparams{\LaTeXparams #1}}%
\newif\ifBoxFrame \BoxFramefalse
\newif\ifOverFrame \OverFramefalse
\newif\ifUnderFrame \UnderFramefalse

\def\BOXTHEFRAME#1{%
   \hbox{%
      \ifBoxFrame
         \frame{#1}%
      \else
         {#1}%
      \fi
   }%
}

\def\doFRAMEparams#1{\BoxFramefalse\OverFramefalse\UnderFramefalse\readFRAMEparams#1\end}%
\def\readFRAMEparams#1{%
 \ifx#1\end%
  \let\next=\relax
  \else
  \ifx#1i\dispkind=\z@\fi
  \ifx#1d\dispkind=\@ne\fi
  \ifx#1f\dispkind=\tw@\fi
  \ifx#1t\addtoLaTeXparams{t}\fi
  \ifx#1b\addtoLaTeXparams{b}\fi
  \ifx#1p\addtoLaTeXparams{p}\fi
  \ifx#1h\addtoLaTeXparams{h}\fi
  \ifx#1X\BoxFrametrue\fi
  \ifx#1O\OverFrametrue\fi
  \ifx#1U\UnderFrametrue\fi
  \ifx#1w
    \ifnum\draft=1\wasdrafttrue\else\wasdraftfalse\fi
    \draft=\@ne
  \fi
  \let\next=\readFRAMEparams
  \fi
 \next
 }%
\def\IFRAME#1#2#3#4#5#6{%
      \bgroup
      \let\QCTOptA\empty
      \let\QCTOptB\empty
      \let\QCBOptA\empty
      \let\QCBOptB\empty
      #6%
      \parindent=0pt%
      \leftskip=0pt
      \rightskip=0pt
      \setbox0 = \hbox{\QCBOptA}%
      \@tempdima = #1\relax
      \ifOverFrame
          \typeout{This is not implemented yet}%
          \show\HELP
      \else
         \ifdim\wd0>\@tempdima
            \advance\@tempdima by \@tempdima
            \ifdim\wd0 >\@tempdima
               \textwidth=\@tempdima
               \setbox1 =\vbox{%
                  \noindent\hbox to \@tempdima{\hfill\GRAPHIC{#5}{#4}{#1}{#2}{#3}\hfill}\\%
                  \noindent\hbox to \@tempdima{\parbox[b]{\@tempdima}{\QCBOptA}}%
               }%
               \wd1=\@tempdima
            \else
               \textwidth=\wd0
               \setbox1 =\vbox{%
                 \noindent\hbox to \wd0{\hfill\GRAPHIC{#5}{#4}{#1}{#2}{#3}\hfill}\\%
                 \noindent\hbox{\QCBOptA}%
               }%
               \wd1=\wd0
            \fi
         \else
            \ifdim\wd0>0pt
              \hsize=\@tempdima
              \setbox1 =\vbox{%
                \unskip\GRAPHIC{#5}{#4}{#1}{#2}{0pt}%
                \break
                \unskip\hbox to \@tempdima{\hfill \QCBOptA\hfill}%
              }%
              \wd1=\@tempdima
           \else
              \hsize=\@tempdima
              \setbox1 =\vbox{%
                \unskip\GRAPHIC{#5}{#4}{#1}{#2}{0pt}%
              }%
              \wd1=\@tempdima
           \fi
         \fi
         \@tempdimb=\ht1
         \advance\@tempdimb by \dp1
         \advance\@tempdimb by -#2%
         \advance\@tempdimb by #3%
         \leavevmode
         \raise -\@tempdimb \hbox{\box1}%
      \fi
      \egroup%
}%
\def\DFRAME#1#2#3#4#5{%
 \begin{center}
     \let\QCTOptA\empty
     \let\QCTOptB\empty
     \let\QCBOptA\empty
     \let\QCBOptB\empty
     \ifOverFrame 
        #5\QCTOptA\par
     \fi
     \GRAPHIC{#4}{#3}{#1}{#2}{\z@}
     \ifUnderFrame 
        \nobreak\par #5\QCBOptA
     \fi
 \end{center}%
 }%
\def\FFRAME#1#2#3#4#5#6#7{%
 \begin{figure}[#1]%
  \let\QCTOptA\empty
  \let\QCTOptB\empty
  \let\QCBOptA\empty
  \let\QCBOptB\empty
  \ifOverFrame
    #4
    \ifx\QCTOptA\empty
    \else
      \ifx\QCTOptB\empty
        \caption{\QCTOptA}%
      \else
        \caption[\QCTOptB]{\QCTOptA}%
      \fi
    \fi
    \ifUnderFrame\else
      \label{#5}%
    \fi
  \else
    \UnderFrametrue%
  \fi
  \begin{center}\GRAPHIC{#7}{#6}{#2}{#3}{\z@}\end{center}%
  \ifUnderFrame
    #4
    \ifx\QCBOptA\empty
      \caption{}%
    \else
      \ifx\QCBOptB\empty
        \caption{\QCBOptA}%
      \else
        \caption[\QCBOptB]{\QCBOptA}%
      \fi
    \fi
    \label{#5}%
  \fi
  \end{figure}%
 }%
\newcount\dispkind%

\def\makeactives{
  \catcode`\"=\active
  \catcode`\;=\active
  \catcode`\:=\active
  \catcode`\'=\active
  \catcode`\~=\active
}
\bgroup
   \makeactives
   \gdef\activesoff{%
      \def"{\string"}
      \def;{\string;}
      \def:{\string:}
      \def'{\string'}
      \def~{\string~}
    }
\egroup

\def\FRAME#1#2#3#4#5#6#7#8{%
 \bgroup
 \@ifundefined{bbl@deactivate}{}{\activesoff}
 \ifnum\draft=\@ne
   \wasdrafttrue
 \else
   \wasdraftfalse%
 \fi
 \def\LaTeXparams{}%
 \dispkind=\z@
 \def\LaTeXparams{}%
 \doFRAMEparams{#1}%
 \ifnum\dispkind=\z@\IFRAME{#2}{#3}{#4}{#7}{#8}{#5}\else
  \ifnum\dispkind=\@ne\DFRAME{#2}{#3}{#7}{#8}{#5}\else
   \ifnum\dispkind=\tw@
    \edef\@tempa{\noexpand\FFRAME{\LaTeXparams}}%
    \@tempa{#2}{#3}{#5}{#6}{#7}{#8}%
    \fi
   \fi
  \fi
  \ifwasdraft\draft=1\else\draft=0\fi{}%
  \egroup
 }%
\def\TEXUX#1{"texux"}

\long\def\QQQ#1#2{%
     \long\expandafter\def\csname#1\endcsname{#2}}%
\@ifundefined{QTP}{\def\QTP#1{}}{}
\@ifundefined{QEXCLUDE}{\def\QEXCLUDE#1{}}{}
\@ifundefined{Qlb}{}{}
\@ifundefined{Qlt}{}{}
\long\def\QQA#1#2{}%
\def\QTR#1#2{{\csname#1\endcsname #2}}
\def\EXPAND#1[#2]#3{}%
\def\NOEXPAND#1[#2]#3{}%
\def\LaTeXparent#1{}%
\def\ChildStyles#1{}%
\def\ChildDefaults#1{}%
\def\QTagDef#1#2#3{}%
\@ifundefined{StyleEditBeginDoc}{}{}
\def\QQfnmark#1{\footnotemark}

\@ifundefined{INDEX}{\def\INDEX#1#2{}{}}{}%
\@ifundefined{SUBINDEX}{\def\SUBINDEX#1#2#3{}{}{}}{}%
\@ifundefined{initial}%
   {\def\initial#1{\bigbreak{\raggedright\large\bf #1}\kern 2\p@\penalty3000}}%
   {}%
\@ifundefined{entry}{}{}%
\@ifundefined{primary}{}{}%
\@ifundefined{secondary}{}{}%
\@ifundefined{ZZZ}{}{\makeatletter\input gnuindex.sty\makeatother\makeindex\makeatletter}%
\@ifundefined{abstract}{%
 \def\abstract{%
  \if@twocolumn
   \section*{Abstract (Not appropriate in this style!)}%
   \else \small 
   \begin{center}{\bf Abstract\vspace{-.5em}\vspace{\z@}}\end{center}%
   \quotation 
   \fi
  }%
 }{%
 }%
\@ifundefined{endabstract}{\def\endabstract
  {\if@twocolumn\else\endquotation\fi}}{}%
\@ifundefined{maketitle}{\def\maketitle#1{}}{}%
\@ifundefined{affiliation}{\def\affiliation#1{}}{}%
\@ifundefined{proof}{}{}%
\@ifundefined{endproof}{}{}%
\@ifundefined{newfield}{\def\newfield#1#2{}}{}%
\@ifundefined{chapter}{\def\chapter#1{\par(Chapter head:)#1\par }%
 \newcount\c@chapter}{}%
\@ifundefined{part}{\def\part#1{\par(Part head:)#1\par }}{}%
\@ifundefined{section}{\def\section#1{\par(Section head:)#1\par }}{}%
\@ifundefined{subsection}{\def\subsection#1%
 {\par(Subsection head:)#1\par }}{}%
\@ifundefined{subsubsection}{\def\subsubsection#1%
 {\par(Subsubsection head:)#1\par }}{}%
\@ifundefined{paragraph}{\def\paragraph#1%
 {\par(Subsubsubsection head:)#1\par }}{}%
\@ifundefined{subparagraph}{\def\subparagraph#1%
 {\par(Subsubsubsubsection head:)#1\par }}{}%
\@ifundefined{therefore}{}{}%
\@ifundefined{backepsilon}{}{}%
\@ifundefined{yen}{}{}%
\@ifundefined{registered}{%
   \def\registered{\relax\ifmmode{}\r@gistered
                    \else$\m@th\r@gistered$\fi}%
 \def\r@gistered{^{\ooalign
  {\hfil\raise.07ex\hbox{$\scriptstyle\rm\text{R}$}\hfil\crcr
  \mathhexbox20D}}}}{}%
\@ifundefined{Eth}{}{}%
\@ifundefined{eth}{}{}%
\@ifundefined{Thorn}{}{}%
\@ifundefined{thorn}{}{}%
\@ifundefined{degree}{}{}%
\newdimen\theight
\def\Column{%
 \vadjust{\setbox\z@=\hbox{\scriptsize\quad\quad tcol}%
  \theight=\ht\z@\advance\theight by \dp\z@\advance\theight by \lineskip
  \kern -\theight \vbox to \theight{%
   \rightline{\rlap{\box\z@}}%
   \vss
   }%
  }%
 }%
\def\qed{%
 \ifhmode\unskip\nobreak\fi\ifmmode\ifinner\else\hskip5\p@\fi\fi
 \hbox{\hskip5\p@\vrule width4\p@ height6\p@ depth1.5\p@\hskip\p@}%
 }%
\def\miss{\hbox{\vrule height2\p@ width 2\p@ depth\z@}}%
\def\tcol#1{{\baselineskip=6\p@ \vcenter{#1}} \Column}  %
\def\newfmtname{LaTeX2e}
\def\chkcompat{%
   \if@compatibility
   \else
     \usepackage{latexsym}
   \fi
}

\ifx\fmtname\newfmtname
  \DeclareOldFontCommand{\rm}{\normalfont\rmfamily}{\mathrm}
  \DeclareOldFontCommand{\sf}{\normalfont\sffamily}{\mathsf}
  \DeclareOldFontCommand{\tt}{\normalfont\ttfamily}{\mathtt}
  \DeclareOldFontCommand{\bf}{\normalfont\bfseries}{\mathbf}
  \DeclareOldFontCommand{\it}{\normalfont\itshape}{\mathit}
  \DeclareOldFontCommand{\sl}{\normalfont\slshape}{\@nomath\sl}
  \DeclareOldFontCommand{\sc}{\normalfont\scshape}{\@nomath\sc}
  \chkcompat
\fi

\def\alpha{{\Greekmath 010B}}%
\def\beta{{\Greekmath 010C}}%
\def\gamma{{\Greekmath 010D}}%
\def\delta{{\Greekmath 010E}}%
\def\epsilon{{\Greekmath 010F}}%
\def\zeta{{\Greekmath 0110}}%
\def\eta{{\Greekmath 0111}}%
\def\theta{{\Greekmath 0112}}%
\def\iota{{\Greekmath 0113}}%
\def\kappa{{\Greekmath 0114}}%
\def\lambda{{\Greekmath 0115}}%
\def\mu{{\Greekmath 0116}}%
\def\nu{{\Greekmath 0117}}%
\def\xi{{\Greekmath 0118}}%
\def\pi{{\Greekmath 0119}}%
\def\rho{{\Greekmath 011A}}%
\def\sigma{{\Greekmath 011B}}%
\def\tau{{\Greekmath 011C}}%
\def\upsilon{{\Greekmath 011D}}%
\def\phi{{\Greekmath 011E}}%
\def\chi{{\Greekmath 011F}}%
\def\psi{{\Greekmath 0120}}%
\def\omega{{\Greekmath 0121}}%
\def\varepsilon{{\Greekmath 0122}}%
\def\vartheta{{\Greekmath 0123}}%
\def\varpi{{\Greekmath 0124}}%
\def\varrho{{\Greekmath 0125}}%
\def\varsigma{{\Greekmath 0126}}%
\def\varphi{{\Greekmath 0127}}%

\def\nabla{{\Greekmath 0272}}
\def\FindBoldGroup{%
   {\setbox0=\hbox{$\mathbf{x\global\edef\theboldgroup{\the\mathgroup}}$}}%
}

\def\Greekmath#1#2#3#4{%
    \if@compatibility
        \ifnum\mathgroup=\symbold
           \mathchoice{\mbox{\boldmath$\displaystyle\mathchar"#1#2#3#4$}}%
                      {\mbox{\boldmath$\textstyle\mathchar"#1#2#3#4$}}%
                      {\mbox{\boldmath$\scriptstyle\mathchar"#1#2#3#4$}}%
                      {\mbox{\boldmath$\scriptscriptstyle\mathchar"#1#2#3#4$}}%
        \else
           \mathchar"#1#2#3#4%
        \fi 
    \else 
        \FindBoldGroup
        \ifnum\mathgroup=\theboldgroup 
           \mathchoice{\mbox{\boldmath$\displaystyle\mathchar"#1#2#3#4$}}%
                      {\mbox{\boldmath$\textstyle\mathchar"#1#2#3#4$}}%
                      {\mbox{\boldmath$\scriptstyle\mathchar"#1#2#3#4$}}%
                      {\mbox{\boldmath$\scriptscriptstyle\mathchar"#1#2#3#4$}}%
        \else
           \mathchar"#1#2#3#4%
        \fi     	    
	  \fi}

\newif\ifGreekBold  \GreekBoldfalse
\let\SAVEPBF=\pbf
\def\pbf{\GreekBoldtrue\SAVEPBF}%

\@ifundefined{theorem}{}{}
\@ifundefined{lemma}{}{}
\@ifundefined{corollary}{}{}
\@ifundefined{conjecture}{}{}
\@ifundefined{proposition}{}{}
\@ifundefined{axiom}{}{}
\@ifundefined{remark}{\newtheorem{remark}{Remark}}{}
\@ifundefined{example}{\newtheorem{example}{Example}}{}
\@ifundefined{exercise}{}{}
\@ifundefined{definition}{}{}

\@ifundefined{mathletters}{%
  \newcounter{equationnumber}  
  \def\mathletters{%
     \addtocounter{equation}{1}
     \edef\@currentlabel{\theequation}%
     \setcounter{equationnumber}{\c@equation}
     \setcounter{equation}{0}%
     \edef\theequation{\@currentlabel\noexpand\alph{equation}}%
  }
  
}{}

\@ifundefined{BibTeX}{%
    \def\BibTeX{{\rm B\kern-.05em{\sc i\kern-.025em b}\kern-.08em
                 T\kern-.1667em\lower.7ex\hbox{E}\kern-.125emX}}}{}%
\@ifundefined{AmS}%
    {\def\AmS{{\protect\usefont{OMS}{cmsy}{m}{n}%
                A\kern-.1667em\lower.5ex\hbox{M}\kern-.125emS}}}{}%
\@ifundefined{AmSTeX}{}{}%
\ifx\ds@amstex\relax
\makeatother\endinput
\else
   \@ifpackageloaded{amstex}%
      {\message{amstex already loaded}\makeatother\endinput}
      {}
   \@ifpackageloaded{amsgen}%
      {\message{amsgen already loaded}\makeatother\endinput}
      {}
\fi
\let\DOTSI\relax
\def\RIfM@{\relax\ifmmode}%
\def\FN@{\futurelet\next}%
\newcount\intno@
\def\iint{\DOTSI\intno@\tw@\FN@\ints@}%
\def\iiint{\DOTSI\intno@\thr@@\FN@\ints@}%
\def\iiiint{\DOTSI\intno@4 \FN@\ints@}%
\def\idotsint{\DOTSI\intno@\z@\FN@\ints@}%
\def\ints@{\findlimits@\ints@@}%
\newif\iflimtoken@
\newif\iflimits@
\def\findlimits@{\limtoken@true\ifx\next\limits\limits@true
 \else\ifx\next\nolimits\limits@false\else
 \limtoken@false\ifx\ilimits@\nolimits\limits@false\else
 \ifinner\limits@false\else\limits@true\fi\fi\fi\fi}%
\def\multint@{\int\ifnum\intno@=\z@\intdots@                          
 \else\intkern@\fi                                                    
 \ifnum\intno@>\tw@\int\intkern@\fi                                   
 \ifnum\intno@>\thr@@\int\intkern@\fi                                 
 \int}
\def\multintlimits@{\intop\ifnum\intno@=\z@\intdots@\else\intkern@\fi
 \ifnum\intno@>\tw@\intop\intkern@\fi
 \ifnum\intno@>\thr@@\intop\intkern@\fi\intop}%
\def\intic@{%
    \mathchoice{\hskip.5em}{\hskip.4em}{\hskip.4em}{\hskip.4em}}%
\def\negintic@{\mathchoice
 {\hskip-.5em}{\hskip-.4em}{\hskip-.4em}{\hskip-.4em}}%
\def\ints@@{\iflimtoken@                                              
 \def\ints@@@{\iflimits@\negintic@
   \mathop{\intic@\multintlimits@}\limits                             
  \else\multint@\nolimits\fi                                          
  \eat@}
 \else                                                                
 \def\ints@@@{\iflimits@\negintic@
  \mathop{\intic@\multintlimits@}\limits\else
  \multint@\nolimits\fi}\fi\ints@@@}%
\def\intkern@{\mathchoice{\!\!\!}{\!\!}{\!\!}{\!\!}}%
\def\plaincdots@{\mathinner{\cdotp\cdotp\cdotp}}%
\def\intdots@{\mathchoice{\plaincdots@}%
 {{\cdotp}\mkern1.5mu{\cdotp}\mkern1.5mu{\cdotp}}%
 {{\cdotp}\mkern1mu{\cdotp}\mkern1mu{\cdotp}}%
 {{\cdotp}\mkern1mu{\cdotp}\mkern1mu{\cdotp}}}%
\def\RIfM@{\relax\protect\ifmmode}
\def\text{\RIfM@\expandafter\text@\else\expandafter\mbox\fi}
\let\nfss@text\text
\def\text@#1{\mathchoice
   {\textdef@\displaystyle\f@size{#1}}%
   {\textdef@\textstyle\tf@size{\firstchoice@false #1}}%
   {\textdef@\textstyle\sf@size{\firstchoice@false #1}}%
   {\textdef@\textstyle \ssf@size{\firstchoice@false #1}}%
   \glb@settings}

\def\textdef@#1#2#3{\hbox{{%
                    \everymath{#1}%
                    \let\f@size#2\selectfont
                    #3}}}
\newif\iffirstchoice@
\firstchoice@true
\def\Let@{\relax\iffalse{\fi\let\\=\cr\iffalse}\fi}%
\def\vspace@{\def\vspace##1{\crcr\noalign{\vskip##1\relax}}}%
\def\multilimits@{\bgroup\vspace@\Let@
 \baselineskip\fontdimen10 \scriptfont\tw@
 \advance\baselineskip\fontdimen12 \scriptfont\tw@
 \lineskip\thr@@\fontdimen8 \scriptfont\thr@@
 \lineskiplimit\lineskip
 \vbox\bgroup\ialign\bgroup\hfil$\m@th\scriptstyle{##}$\hfil\crcr}%
\def\Sb{_\multilimits@}%
\def\endSb{\crcr\egroup\egroup\egroup}%
\def\Sp{^\multilimits@}%

\newdimen\ex@
\def\rightarrowfill@#1{$#1\m@th\mathord-\mkern-6mu\cleaders
 \hbox{$#1\mkern-2mu\mathord-\mkern-2mu$}\hfill
 \mkern-6mu\mathord\rightarrow$}%
\def\leftarrowfill@#1{$#1\m@th\mathord\leftarrow\mkern-6mu\cleaders
 \hbox{$#1\mkern-2mu\mathord-\mkern-2mu$}\hfill\mkern-6mu\mathord-$}%
\def\leftrightarrowfill@#1{$#1\m@th\mathord\leftarrow
\mkern-6mu\cleaders
 \hbox{$#1\mkern-2mu\mathord-\mkern-2mu$}\hfill
 \mkern-6mu\mathord\rightarrow$}%
\def\overrightarrow{\mathpalette\overrightarrow@}%
\def\overrightarrow@#1#2{\vbox{\ialign{##\crcr\rightarrowfill@#1\crcr
 \noalign{\kern-\ex@\nointerlineskip}$\m@th\hfil#1#2\hfil$\crcr}}}%

\def\overleftarrow{\mathpalette\overleftarrow@}%
\def\overleftarrow@#1#2{\vbox{\ialign{##\crcr\leftarrowfill@#1\crcr
 \noalign{\kern-\ex@\nointerlineskip}$\m@th\hfil#1#2\hfil$\crcr}}}%
\def\overleftrightarrow{\mathpalette\overleftrightarrow@}%
\def\overleftrightarrow@#1#2{\vbox{\ialign{##\crcr
   \leftrightarrowfill@#1\crcr
 \noalign{\kern-\ex@\nointerlineskip}$\m@th\hfil#1#2\hfil$\crcr}}}%
\def\underrightarrow{\mathpalette\underrightarrow@}%
\def\underrightarrow@#1#2{\vtop{\ialign{##\crcr$\m@th\hfil#1#2\hfil
  $\crcr\noalign{\nointerlineskip}\rightarrowfill@#1\crcr}}}%

\def\underleftarrow{\mathpalette\underleftarrow@}%
\def\underleftarrow@#1#2{\vtop{\ialign{##\crcr$\m@th\hfil#1#2\hfil
  $\crcr\noalign{\nointerlineskip}\leftarrowfill@#1\crcr}}}%
\def\underleftrightarrow{\mathpalette\underleftrightarrow@}%
\def\underleftrightarrow@#1#2{\vtop{\ialign{##\crcr$\m@th
  \hfil#1#2\hfil$\crcr
 \noalign{\nointerlineskip}\leftrightarrowfill@#1\crcr}}}%

\def\qopnamewl@#1{\mathop{\operator@font#1}\nlimits@}
\let\nlimits@\displaylimits
\def\setboxz@h{\setbox\z@\hbox}

\def\varlim@#1#2{\mathop{\vtop{\ialign{##\crcr
 \hfil$#1\m@th\operator@font lim$\hfil\crcr
 \noalign{\nointerlineskip}#2#1\crcr
 \noalign{\nointerlineskip\kern-\ex@}\crcr}}}}

 \def\rightarrowfill@#1{\m@th\setboxz@h{$#1-$}\ht\z@\z@
  $#1\copy\z@\mkern-6mu\cleaders
  \hbox{$#1\mkern-2mu\box\z@\mkern-2mu$}\hfill
  \mkern-6mu\mathord\rightarrow$}
\def\leftarrowfill@#1{\m@th\setboxz@h{$#1-$}\ht\z@\z@
  $#1\mathord\leftarrow\mkern-6mu\cleaders
  \hbox{$#1\mkern-2mu\copy\z@\mkern-2mu$}\hfill
  \mkern-6mu\box\z@$}

\def\projlim{\qopnamewl@{proj\,lim}}
\def\injlim{\qopnamewl@{inj\,lim}}
\def\varinjlim{\mathpalette\varlim@\rightarrowfill@}
\def\varprojlim{\mathpalette\varlim@\leftarrowfill@}
\def\varliminf{\mathpalette\varliminf@{}}
\def\varliminf@#1{\mathop{\underline{\vrule\@depth.2\ex@\@width\z@
   \hbox{$#1\m@th\operator@font lim$}}}}
\def\varlimsup{\mathpalette\varlimsup@{}}
\def\varlimsup@#1{\mathop{\overline
  {\hbox{$#1\m@th\operator@font lim$}}}}

\begingroup \catcode `|=0 \catcode `[= 1
\catcode`]=2 \catcode `\{=12 \catcode `\}=12
\catcode`\\=12 
|gdef|@alignverbatim#1\end{align}[#1|end[align]]
|gdef|@salignverbatim#1\end{align*}[#1|end[align*]]

|gdef|@alignatverbatim#1\end{alignat}[#1|end[alignat]]
|gdef|@salignatverbatim#1\end{alignat*}[#1|end[alignat*]]

|gdef|@xalignatverbatim#1\end{xalignat}[#1|end[xalignat]]
|gdef|@sxalignatverbatim#1\end{xalignat*}[#1|end[xalignat*]]

|gdef|@gatherverbatim#1\end{gather}[#1|end[gather]]
|gdef|@sgatherverbatim#1\end{gather*}[#1|end[gather*]]

|gdef|@gatherverbatim#1\end{gather}[#1|end[gather]]
|gdef|@sgatherverbatim#1\end{gather*}[#1|end[gather*]]

|gdef|@multilineverbatim#1\end{multiline}[#1|end[multiline]]
|gdef|@smultilineverbatim#1\end{multiline*}[#1|end[multiline*]]

|gdef|@arraxverbatim#1\end{arrax}[#1|end[arrax]]
|gdef|@sarraxverbatim#1\end{arrax*}[#1|end[arrax*]]

|gdef|@tabulaxverbatim#1\end{tabulax}[#1|end[tabulax]]
|gdef|@stabulaxverbatim#1\end{tabulax*}[#1|end[tabulax*]]

|endgroup

\def\align{\@verbatim \frenchspacing\@vobeyspaces \@alignverbatim
You are using the "align" environment in a style in which it is not defined.}

\@namedef{align*}{\@verbatim\@salignverbatim
You are using the "align*" environment in a style in which it is not defined.}
\expandafter\let\csname endalign*\endcsname =\endtrivlist

\def\alignat{\@verbatim \frenchspacing\@vobeyspaces \@alignatverbatim
You are using the "alignat" environment in a style in which it is not defined.}

\@namedef{alignat*}{\@verbatim\@salignatverbatim
You are using the "alignat*" environment in a style in which it is not defined.}
\expandafter\let\csname endalignat*\endcsname =\endtrivlist

\def\xalignat{\@verbatim \frenchspacing\@vobeyspaces \@xalignatverbatim
You are using the "xalignat" environment in a style in which it is not defined.}

\@namedef{xalignat*}{\@verbatim\@sxalignatverbatim
You are using the "xalignat*" environment in a style in which it is not defined.}
\expandafter\let\csname endxalignat*\endcsname =\endtrivlist

\def\gather{\@verbatim \frenchspacing\@vobeyspaces \@gatherverbatim
You are using the "gather" environment in a style in which it is not defined.}

\@namedef{gather*}{\@verbatim\@sgatherverbatim
You are using the "gather*" environment in a style in which it is not defined.}
\expandafter\let\csname endgather*\endcsname =\endtrivlist

\def\multiline{\@verbatim \frenchspacing\@vobeyspaces \@multilineverbatim
You are using the "multiline" environment in a style in which it is not defined.}

\@namedef{multiline*}{\@verbatim\@smultilineverbatim
You are using the "multiline*" environment in a style in which it is not defined.}
\expandafter\let\csname endmultiline*\endcsname =\endtrivlist

\def\arrax{\@verbatim \frenchspacing\@vobeyspaces \@arraxverbatim
You are using a type of "array" construct that is only allowed in AmS-LaTeX.}

\def\tabulax{\@verbatim \frenchspacing\@vobeyspaces \@tabulaxverbatim
You are using a type of "tabular" construct that is only allowed in AmS-LaTeX.}

\@namedef{arrax*}{\@verbatim\@sarraxverbatim
You are using a type of "array*" construct that is only allowed in AmS-LaTeX.}
\expandafter\let\csname endarrax*\endcsname =\endtrivlist

\@namedef{tabulax*}{\@verbatim\@stabulaxverbatim
You are using a type of "tabular*" construct that is only allowed in AmS-LaTeX.}
\expandafter\let\csname endtabulax*\endcsname =\endtrivlist

\def\@@eqncr{\let\@tempa\relax
    \ifcase\@eqcnt \def\@tempa{& & &}\or \def\@tempa{& &}%
      \else \def\@tempa{&}\fi
     \@tempa
     \if@eqnsw
        \iftag@
           \@taggnum
        \else
           \@eqnnum\stepcounter{equation}%
        \fi
     \fi
     \global\tag@false
     \global\@eqnswtrue
     \global\@eqcnt\z@\cr}

 \def\endequation{%
     \ifmmode\ifinner 
      \iftag@
        \addtocounter{equation}{-1} 
        $\hfil
           \displaywidth\linewidth\@taggnum\egroup \endtrivlist
        \global\tag@false
        \global\@ignoretrue   
      \else
        $\hfil
           \displaywidth\linewidth\@eqnnum\egroup \endtrivlist
        \global\tag@false
        \global\@ignoretrue 
      \fi
     \else   
      \iftag@
        \addtocounter{equation}{-1} 
        \eqno \hbox{\@taggnum}
        \global\tag@false%
        $$\global\@ignoretrue
      \else
        \eqno \hbox{\@eqnnum}
        $$\global\@ignoretrue
      \fi
     \fi\fi
 } 

 \newif\iftag@ \tag@false
 
 \def\tag{\@ifnextchar*{\@tagstar}{\@tag}}
 \def\@tag#1{%
     \global\tag@true
     \global\def\@taggnum{(#1)}}
 \def\@tagstar*#1{%
     \global\tag@true
     \global\def\@taggnum{#1}%
}

\makeatother

\makeatletter
\input{tcilatex}
\input tcilatex
\makeatother
\begin{document}

\author{A. I. Fedoseyev$^{*}$, M. J. Friedman$^{\dagger }$ and E. J. Kansa$%
^{\ddagger }$ \\
$^{*}$Center for Microgravity and Materials Research, \\
University of Alabama in Huntsville, Huntsville, AL 35899\\
$^{\dagger }$Department of Mathematical Sciences, \\
University of Alabama in Huntsville, Huntsville, AL 35899\\
$^{\ddagger }$Lawrence Livermore National Laboratory, Livermore, CA 94551}
\title{Continuation for Nonlinear Elliptic Partial Differential Equations
Discretized by the Multiquadric Method }
\maketitle

\begin{abstract}
The Multiquadric Radial Basis Function (MQ) Method is a meshless collocation
method with global basis functions. It is known to have exponentional
convergence for interpolation problems. We descretize nonlinear elliptic
PDEs by the MQ method. This results in modest size systems of nonlinear
algebraic equations which can be efficiently continued by standard
continuation software such as \textsc{auto} and \textsc{content}. Examples
are given of detection of bifurcations in 1D and 2D PDEs. These examples
show high accuracy with small number of unknowns, as compared with known
results from the literature.

\end{abstract}

\section{Introduction}

Nonlinear multidimensional elliptic partial differential equations (PDEs)
are the basis for many scientific and engineering problems, such as pattern
formation in biology, viscous fluid flow phenomena, chemical reactions,
crystal growth processes, etc. In these problems it is crucial to understand
the qualitative dependence of the solution on the problem parameters. 

During the past two decades the numerical continuation approach has become
popular for qualitative study of solutions to nonlinear equations, see e.g. 
\cite{Rheinboldt86}, \cite{DoKelKer91}, \cite{Seydel98} and references
therein. Several software packages, such as \textsc{auto97} \cite{DCFKSWF97}
and \textsc{content} \cite{KuzLev98}, are currently available for
bifurcation analysis of systems of nonlinear algebraic equations and ODEs,
with only limited bifurcation analysis for 1D elliptic PDEs. For 2D PDEs, we
mention the software package \textsc{pltmg} \cite{Ba:98} that allows to
solve a class of boundary value problems on regions in the plane, to
continue the solution with respect to a parameter and even to compute limit
and branching points. This software combines a sophisticated finite element
discretization with advanced linear algebra techniques. Numerical
continuation for 1D and 2D elliptic PDEs is currently an active research
area, see e.g. \cite{Ne:93}, \cite{ShKe:93}, \cite{SchTimLos96}, \cite
{ChShMe97}, \cite{Davidson97}, \cite{KSLGS98}, \cite{ChiChe98}, and 
\cite[Ch 10]{Go99} for reaction diffusion equations; and \cite{Polia95}, 
\cite{ChTu95} for CFD. The typical approaches used are based on the finite
element or finite difference discretization of the PDEs. They result in very
large (thousands or tens of thousands for 2D problems) systems of nonlinear
algebraic equations with sparse matrices. The continuation process is
typically based on the predictor-corrector algorithms that require solving
nonlinear systems by the Newton type method at each continuation step. For
the bifurcation analysis during the continuation process, one usually needs
to compute at least few eigenvalues of the Jacobian matrix at each
continuation step. The methods currently used both for the continuation and
the corresponding eigenvalue problems are variants of Krylov subspace
methods and recursive projection (RPM). Solving the resulting linear system
and the eigenvalue problem require sophisticated algorithms and considerable
computer resources (CPU time, memory, disk space, etc.).

In this paper we report results of numerical experiments with continuation
and detection of bifurcations for 1D and 2D elliptic PDEs discretized by the
Multiquadric Radial Basis Function (MQ) method. The MQ method was first
introduced for solving PDEs in 1990 by Kansa \cite{Kansa90a}, \cite{Kansa90b}%
. It is a meshless collocation method with global basis functions which
leads to finite dimensional problems with full matrices. It was shown to
give very high accuracy with a relatively small number of unknowns (tens or
hundreds for 2D problems). The corresponding linear systems can be
efficiently solved by direct methods. This opens a possibility for using
standard continuation software, such as \textsc{auto} and \textsc{content},
designed for bifurcation analysis of modest size problems. We also note that
the MQ method does not require predetermined location of the nodes as the
spectral method does.

In Section \ref{review} we summarize previous results on solving PDEs by the
MQ method and our experiments with an eigenvalue problem.

In Section \ref{MQ-method} we formulate an adaptation of the MQ\ method for
the discretization of the parametrized elliptic PDEs.

In Section \ref{examples} we present results of our numerical experiments
with continuation of solutions and detection of bifurcations for 1D and 2D
elliptic PDEs.

In Section \ref{conclusions} we discuss our results.

\section{Review of multiquadric method for elliptic PDEs}

\label{review}

\subsection{Summary of previous results}

The concept of solving PDEs using the radial basis functions (RBF) was
introduced by Kansa in 1990 \cite{Kansa90a}, \cite{Kansa90b}. He implemented
this approach for the solution of hyperbolic, parabolic, and elliptic PDEs
using the MQ RBFs proposed by Hardy \cite{Hardy71}, \cite{Hardy90} for
interpolation of scattered data.

There exists an infinite class of RBFs. A radial basis function, $f(x)$, $%
x\in \mathbb{R}^{n}$, depends only upon the distances between the nodes. A
MQ RBF is $g_{j}(x)=((x-x_{j})^{2}+c_{j}^{2})^{1/2}$, where $x_{j}$ is a
reference node and $c_{j}$ is a shape parameter. In the comprehensive study
by Franke \cite{Franke82}, it is shown that MQ RBFs have the excellent
properties for the interpolation of 2D scattered data. Among studied RBFs
still only the MQ RBFs are proven to have the exponential convergence for
the function interpolation \cite{MaNel90}, \cite{WuSha93}.

The numerical experiments by Kansa \cite{Kansa90a}, \cite{Kansa90b}, and
Golberg and Chen \cite{GoCh96} show high efficiency and very accurate
solution with the MQ scheme. Kansa \cite{Kansa90b} showed that MQ method
yields a high accuracy for parabolic and elliptic PDEs. Example for the
transient convection-diffusion problem with steep initial front demonstrated
highly accurate solution by the MQ method with a small number of nodes even
for large cell Peclet number $Pe$. Test cases with $20$ nodes for the MQ
method ran for diffusion coefficient $D$ in the range from $10^{-1}$ to $%
10^{-3}$. The corresponding cell $Pe$ number was from $0.5$ to $50.0$. Exact
and MQ solution are indistinguishable graphically (apparent difference less
than $10^{-4}$) for $D=10^{-1}$ and $10^{-2}$, while small deviation ($5\%$)
was observed at $D=10^{-3},Pe=50$. No instability or wiggles was seen.
Finite difference simulation with $K=200$ nodes and optimal combination of
the central and upwind differences for the case $D=10^{-1},$ $Pe=5$ resulted
in the error of $3\%$, which was still several orders less accurate than the
MQ method solution.

In the numerical experiments with modeling the von Neumann blast wave Kansa \cite
{Kansa90b} compared the exact solution and its derivatives with the MQ
solution ($35$ nodes) and with finite difference ones ($50, 500$ and $5000$
nodes). The error in value and gradients of pressure, density and energy was 
$10^{-6}$ or less for the MQ method, and in the range from $10^{-4}$ to $%
10^{-2}$ for the best finite difference result with $5000$ nodes.

Golberg and Chen \cite{GoCh97} showed that the solution of the 3D Poisson
equation could be obtained with only 60 randomly distributed nodes to the
same degree of accuracy as a FEM solution with 71,000 linear elements. 

Sharan, Kansa, and Gupta \cite{ShKG97} showed that the MQ method yields very
accurate solutions for elliptic PDE problems, including the biharmonic
equation, and that the MQ approach is simple to implement on domains with
irregular boundaries. Dubal et al. \cite{Du93} noted many benefits of using
MQ RBFs to solve the initial value problem for a 3D nonlinear equation for
the collision of two black holes. The resulting discrete system has 2000
unknowns and was solved directly.

Buhmann \cite{Buhm95} showed that RBFs and, in particular, MQ RBFs are
useful for constructing prewavelets and wavelets. Wavelets are most
frequently used in time series analysis, but there are results for using
wavelets to solve PDEs \cite{Fassh98}, \cite{Narcow97}. As Buhmann points
out, one can generate true wavelets by an orthonormalization process. The
wavelets are an elegant way to achieve the same results as multi-grid
schemes. The MQ RBFs are attractive for prewavelet construction due to
exceptional rates of convergence and their infinite differentiability.

Franke and Schaback paper \cite{Franke98} provides the first theoretical
foundation for solving PDEs by collocation using the RBF methods.

Kansa and Hon \cite{KaHo98} studied several methods to solve linear
equations that arise from the MQ collocation problems. They studied the 2D
Poisson equation, and showed that ill-conditioning of the system of
equations could be circumvented by using the sub-structuring methods.

Kansa \cite{Kansa90b} introduced the concept of variable shape parameters $%
c_{j}$ in the MQ scheme that appeared to work well in some cases. In the
work by Kansa and Hon \cite{KaHo98}, a recipe based upon the local radius of
curvature of solution surface was found to perform better than a constant
shape parameter MQ scheme. A simple variable shape parameter formula is
based the local radius of curvature. Kansa and Hon\cite{KaHo98} tested the
MQ method for the 2D Poisson equation with a set of exact solutions like $%
F=\exp (ax+by),\cos (ax+by),\sin (ax+by)$, $\log (ax+by+c)$, $\exp
(-a(x-1/2)2-b(y-1/2)2)$ and $\arctan (ax+by)$. They showed obtained a high
accuracy (up to $10^{-5}$) and a small residual norm ($10^{-4}$) on a modest
node size set (121 nodes) while locally adapting the shape parameter $c_{j}$.

Franke \cite{Franke82} compared (global) RBF interpolation schemes against
many popular compactly supported schemes such as finite difference method,
and found that the global RBF schemes were superior on six criteria.

Madych \cite{Ma92} showed theoretically that the MQ interpolation scheme
converges faster as the constant MQ shape parameter becomes progressively
larger.

The multi-zone method of Wong et al. \cite{Wo98} is yet another alternative
method for improving computational efficiency. This method is readily
parallelizable, and the conditioning of the resulting matrices are much
better.

Hon and Mao \cite{HoMa98} showed that an adaptive algorithm that adjusted
the nodes to follow the peak of the shock wave can produce extremely
accurate results in 1D Burgers equation with only 10 nodes, even for
extremely steep shocks with $Re=10^{4}$.

\subsection{A simple eigenvalue problem.}

Accurate approximation of eigenvalue problems is essential for bifurcation
analysis of PDEs. We have not found references in literature on the
MQ-solution of eigenvalue problems. We therefore present here results for an
eigenvalue problem for 1D Laplace operator. For details on the MQ
discretization see Section \ref{MQ-method}. This is a scalar problem
\begin{eqnarray}
-u^{^{\prime \prime }} &=&\lambda u,  \label{igen-probl} \\
u(0) &=&u(1)=0,  \nonumber
\end{eqnarray}
that has the exact solution: 
\[
\left( \lambda _{m},\,U^{m}(x)\right) =\left( (\pi m)^{2},\text{ }\sin (\pi
mx)\right) ,\text{ }m=1,2,...
\]
where $\left( \lambda _{m},\,U^{m}(x)\right) $ is the $m-th$ eigenpair of (%
\ref{igen-probl}). Introduce the mesh $x_{n}=nh,$ $n=0,1,...,N$, $h=1/N,$
and consider the standard second order finite difference (FDM)
discretization of (\ref{igen-probl}): 
\begin{eqnarray}
-\frac{u_{n+1}-2u_{n}+u_{n-1}}{h^{2}}=\lambda u_{n}, &&\text{ \  \ }%
n=1,...N-1,  \label{igen-probl-FD} \\
u_{0}=u_{N}=0. &&  \nonumber
\end{eqnarray}
The corresponding approximate eigenpairs are given by 
\[
\left( \lambda _{m}^{h},U_{m}^{h}\right) =\left( 4N^{2}\sin ^{2}\frac{\pi m}{%
2N},\text{ }\left[ 
\begin{array}{l}
\sin \frac{\pi m}{N} \\ 
\sin \frac{\pi 2m}{N} \\ 
.......... \\ 
\sin \frac{\pi (N-1)m}{N}
\end{array}
\right] \right) ,\text{ \  \ }m=1,...,N-1.
\]
We also solved (\ref{igen-probl}) using the MQ discretization for several
values of the number $K$ of internal nodes. Denote by $\left( \lambda
_{m}^{MQ},U_{m}^{MQ}\right) ,$ \  \ $m=1,...,K$ the corresponding
approximate eigenpairs.

The results of our computations are summarized in Table \ref{tab:eigen}. We
use the notation $\varepsilon _{\lambda }^{MQ},$ $\varepsilon _{\lambda }^{h}
$ for the relative errors in $\lambda _{m}^{MQ}$, $\lambda _{m}^{h}$,
respectively, and the notation $\varepsilon _{U}^{MQ}$ for the $L_{\infty }$%
-norm error in $U_{m}^{MQ}$. For each MQ solution we provide a comparison
with the FDM solution that has a sufficient number of nodes to give the same
accuracy for $\lambda _{1}$ as the MQ method. In Part (a) of the table we
use the uniform node distribution for the MQ method. Part (b) of the table
shows that the accuracy of the MQ method can be significantly improved by
adapting the node distribution: we moved only two nodes adjacent to boundary
to reduce their distance from the boundary to $h_{1}=h/4$ (while the
remaining nodes are distributed uniformly).

One can see that the MQ method can give a highly accurate solution with a
small number of unknowns, $10$ $-$ $100$ times smaller than the number of
unknowns in the FDM for the same accuracy.

\begin{table}[tbp] \centering%
\caption{Eigenvalue problem: comparison of eigenvalues and eigenfunctions}%
\label{tab:eigen} \noindent a) MQ method with uniform node distribution for $%
K=5,$ $7\,$and $\,9.$

\begin{tabular}{|c|c|c|c|c|c|}
\hline
$m$ & $\lambda _{m}$ (exact) & $\lambda _{m}^{MQ},\,K=5$ & Rel. err.$%
\,\varepsilon _{\lambda }^{MQ}$ & Rel.err$\,\varepsilon _{U}^{MQ}$ & Rel.
err.$\,\varepsilon _{\lambda }^{h},\,K=47$ \\ \hline
$1$ & $9.86961$ & $9.86596$ & $3.7\times 10^{-4}$ & $3.7\times 10^{-4}$ & $%
3.7\times 10^{-4}$ \\ \hline
$2$ & $39.4784$ & $39.6492$ & $4.3\times 10^{-3}$ & $5.2\times 10^{-3}$ & $%
1.5\times 10^{-3}$ \\ \hline
\end{tabular}

\medskip

\begin{tabular}{|c|c|c|c|c|c|}
\hline
$m$ & $\lambda _{m}$ (exact) & $\lambda _{m}^{MQ},\,K=7$ & Rel. err.$%
\,\varepsilon _{\lambda }^{MQ}$ & Rel. err.$\,\varepsilon _{U}^{MQ}$ & Rel.
err.\thinspace $\varepsilon _{\lambda }^{h},\,K=76$ \\ \hline
$1$ & $9.86961$ & $9.86821$ & $1.4\times 10^{-4}$ & $9.9\times 10^{-5}$ & $%
1.4\times 10^{-4}$ \\ \hline
$2$ & $39.4784$ & $39.4738$ & $1.2\times 10^{-4}$ & $1.8\times 10^{-4}$ & $%
5.7\times 10^{-4}$ \\ \hline
$3$ & $88.8264$ & $89.3648$ & $6.0\times 10^{-3}$ & $1.1\times 10^{-2}$ & $%
1.3\times 10^{-3}$ \\ \hline
\end{tabular}

\medskip

\begin{tabular}{|c|c|c|c|c|c|}
\hline
$m$ & $\lambda _{m}$ (exact) & $\lambda _{m}^{MQ},\,K=9$ & Rel. err.$%
\,\epsilon _{\lambda }^{MQ}$ & Rel. err.$\,\epsilon _{U}^{MQ}$ & Rel. err.$%
\,\epsilon _{\lambda }^{h},\,K=117$ \\ \hline
$1$ & $9.86961$ & $9.86901$ & $6.0\times 10^{-5}$ & $5.0\times 10^{-5}$ & $%
6.0\times 10^{-5}$ \\ \hline
$2$ & $39.4784$ & $39.4846$ & $1.6\times 10^{-4}$ & $2.1\times 10^{-4}$ & $%
2.4\times 10^{-4}$ \\ \hline
$3$ & $88.8264$ & $89.1667$ & $3.8\times 10^{-3}$ & $7.3\times 10^{-3}$ & $%
5.4\times 10^{-4}$ \\ \hline
$4$ & $157.913$ & $159.689$ & $1.1\times 10^{-2}$ & $2.5\times 10^{-2}$ & $%
9.6\times 10^{-4}$ \\ \hline
\end{tabular}

\medskip

\noindent b) MQ method with nonuniform node distribution for $K=7\,$and$\,9$

\begin{tabular}{|c|c|c|c|c|c|}
\hline
$m$ & $\lambda _{m}$ (exact) & $\lambda _{m}^{MQ},\,K=7$ & Rel.err.$%
\,\varepsilon _{\lambda }^{MQ}$ & Rel.err.$\,\varepsilon _{U}^{MQ}$ & 
Rel.err.$\,\varepsilon _{\lambda }^{h},\,K=3477$ \\ \hline
$1$ & $9.86961$ & $9.86961$ & $6.8\times 10^{-8}$ & $3.0\times 10^{-6}$ & $%
6.8\times 10^{-8}$ \\ \hline
$2$ & $39.4784$ & $39.4782$ & $3.2\times 10^{-6}$ & $3.0\times 10^{-4}$ & $%
2.7\times 10^{-7}$ \\ \hline
$3$ & $88.8264$ & $88.8139$ & $1.4\times 10^{-4}$ & $6.5\times 10^{-4}$ & $%
5.4\times 10^{-4}$ \\ \hline
\end{tabular}

\medskip

\begin{tabular}{|c|c|c|c|c|c|}
\hline
$m$ & $\lambda _{m}$ (exact) & $\lambda _{m}^{MQ},\,K=9$ & Rel. err. $%
\epsilon _{\lambda }^{MQ}$ & Rel.err. $\epsilon _{U}^{MQ}$ & Rel.err. $%
\epsilon _{\lambda }^{h},\,K=950$ \\ \hline
$1$ & $9.86961$ & $9.86960$ & $9.1\times 10^{-7}$ & $2.3\times 10^{-6}$ & $%
9.1\times 10^{-7}$ \\ \hline
$2$ & $39.4784$ & $39.4783$ & $1.4\times 10^{-6}$ & $2.0\times 10^{-5}$ & $%
3.6\times 10^{-6}$ \\ \hline
$3$ & $88.8264$ & $88.8241$ & $2.6\times 10^{-5}$ & $1.8\times 10^{-4}$ & $%
8.2\times 10^{-6}$ \\ \hline
$4$ & $157.913$ & $157.882$ & $1.9\times 10^{-4}$ & $1.8\times 10^{-3}$ & $%
1.5\times 10^{-5}$ \\ \hline
\end{tabular}
\end{table}%

\section{Discretization of nonlinear elliptic PDEs by the MQ method}

\label{MQ-method}

We consider the second order system of $n$ parametrized nonlinear elliptic
partial differential equations 
\begin{equation}
D(\alpha )\Delta u-f(\nabla u,u,x,y,\alpha )=0,\;\;\text{  }\alpha \in \Bbb{R%
}\text{, }u(\cdot ),\text{{}}f(\cdot )\in \Bbb{R}^{n}\text{, }(x,y)\in
\Omega \subset \Bbb{R}^{2},  \label{exact_PDE}
\end{equation}
where $D(\alpha )$ is a positive diagonal $n\times n$ matrix, that
dependents smoothly on $\alpha $, subject to boundary conditions 
\begin{equation}
\left. f^{b}(\frac{\partial u}{\partial n},u,x,y,\alpha )\right| _{\partial
\Omega }=0,\text{{}}f^{b}(\cdot )\in \Bbb{R}^{n}.  \label{exact_BC}
\end{equation}
Here $\alpha $ is a control parameter, and we are interested in studying the
dependence of the solutions to the boundary value problem (\ref{exact_PDE}),
(\ref{exact_BC}) on $\alpha .$

We discretize the continuous problem by the multiquadric radial basis
function (MQ) method \cite{Kansa90a}, \cite{Kansa90b}, \cite{MaNel90} as
follows. Introduce a set $\Theta _{h}$ of nodes ($N$ internal and $N_{b}$ on
the boundary) 
\[
\Theta _{h}=\left\{ (x_{i},y_{i})\mid _{i=1,N}\subset \Omega \text{,  }%
(x_{i},y_{i})\mid _{i=N+1,N+N_{b}}\subset \partial \Omega \right\} 
\]
and look for the approximate solution to (\ref{exact_PDE}), (\ref{exact_BC})
in the form\ {} 
\begin{equation}
u_{h}(x,y){=}a_{0}+\sum_{j=1}^{j=N-1}a_{j}\left(
g_{j}(c_{j},x,y)-g_{N}(c_{N},x,y)\right)
+\sum_{j=N+1}^{j=N+N_{b}}a_{j}\left(
g_{j}(c_{j},x,y)-g_{N}(c_{N},x,y)\right) ,  \label{approx_sol}
\end{equation}
where $a_{j}\in \Bbb{R}^{n}$ are the unknown expansion coefficients and
\[
g_{j}(c_{j},x,y)=\sqrt{(x-x_{j})^{2}+(y-y_{j})^{2}+c_{j}^{2}},\;\;\text{ {}}%
j=1,...,N+N_{b},
\]
are the MQ basis functions, and $c_{j}>0$ is called a \textit{shape parameter%
} \cite{Kansa90b}. We then substitute $u_{h}(x,y)$ into (\ref{exact_PDE}), (%
\ref{exact_BC}) and use collocation at the nodes $\Theta {_{h}}$ to obtain a
finite dimensional system 
\begin{equation}
\varphi _{i}(a,\alpha )\equiv D(\alpha )\Delta u_{h}(x_{i},y_{i})-f(\nabla
u_{h}(x_{i},y_{i}),u_{h}(x_{i},y_{i}),x_{i},y_{i},\alpha )=0,\text{ }\;\;%
\text{ }i=1,...,N,  \label{approx_PDE}
\end{equation}
\begin{equation}
\varphi _{i-N}^{b}(a,\alpha )\equiv f^{b}(\frac{\partial u_{h}(x_{i},y_{i})}{%
\partial n},u_{h}(x_{i},y_{i}),x_{i},y_{i},\alpha )=0,\text{ }\;\;\text{ }%
i=N+1,...,N+N_{b}.  \label{approx_BC}
\end{equation}

We next modify the discretized system to make it more suitable for
continuation and bifurcation analysis. 1) We eliminate $a_{j} $, $%
j=N+1,...,N+N_{b} $, associated with the boundary nodes, so as to minimize
the number of unknowns. 2) We reformulate (\ref{approx_sol}) in terms of
nodal values $u_{i} $ so that to have the correct eigenvalue problem (to
avoid dealing with matrix stencils) for the Jacobian matrix of (\ref
{approx_PDE}) for detecting bifurcations during the continuation process.

This is accomplished as follows. Denote $a^{1}=(a_{0},a_{1},...,a_{N-1})\in 
\Bbb{R}^{n\times N}$, $a^{2}=(a_{N+1},...,a_{N+N_{b}})\in \Bbb{R}^{n\times
N_{b}}$, $\varphi =(\varphi _{1},...,\varphi _{N})$, $\varphi ^{b}=(\varphi
_{1}^{b},...,\varphi _{N_{b}}^{b})$, and rewrite the system (\ref{approx_PDE}%
), (\ref{approx_BC}) as 
\begin{equation}
\varphi (a^{1},a^{2},\alpha )=0,\;\;\text{ }\varphi (\cdot )\in \Bbb{R}%
^{n\times N},  \label{approx_PDE_a}
\end{equation}
\begin{equation}
\varphi ^{b}(a^{1},a^{2},\alpha )=0,\;\;\text{ }\varphi ^{b}(\cdot )\in \Bbb{%
R}^{n\times N_{b}}.  \label{approx_BC_a}
\end{equation}
Assuming that the implicit function theorem is applicable here (which is
usually the case), we solve (\ref{approx_BC_a}) for $a^{2}$ to obtain 
\begin{equation}
a^{2}=\psi (a^{1},\alpha ),\text{ or, in components, }a_{j}=\psi
_{j}(a^{1},\alpha ),\text{{}}\;\;\text{ }j=N+1,...,N+N_{b}.  \label{a2_eq}
\end{equation}
Substituting this into (\ref{approx_PDE_a}) yields
\begin{equation}
\varphi (a^{1},\psi (a^{1},\alpha ),\alpha )=0,\;\;\text{ }\varphi (\cdot
)\in \Bbb{R}^{n\times N}\text{.}  \label{fi_eq}
\end{equation}
We next want to reformulate (\ref{fi_eq}) in terms of the nodal values $%
U=(u_{1},u_{2},...,u_{N}$$)\in \Bbb{R}^{n\times N}$ of the approximate
solution at the internal nodes defined by $u_{i}=u_{h}(x_{i},y_{i})$. To
this end we first eliminate $a^{2}=(a_{N+1},...,a_{N+N_{b}})$ from (\ref
{approx_sol}) by substituting (\ref{a2_eq}) into (\ref{approx_sol}) to
obtain 
\begin{equation}
u_{h}(x,y){=}a_{0}+\sum_{j=1}^{j=N-1}a_{j}\left(
g_{j}(c_{j},x,y)-g_{N}(c_{N},x,y)\right) +\sum_{j=N+1}^{j=N+N_{b}}\psi
_{j}(a^{1},\alpha )\left( g_{j}(c_{j},x,y)-g_{N}(c_{N},x,y)\right) 
\end{equation}
We now define the map $\Gamma :$ $a^{1}\longmapsto U=\Gamma (a^{1})$. For $%
i=0,...,N-1:$%
\begin{equation}
u_{i}=a_{0}+\sum_{j=1}^{j=N-1}\left(
g_{j}(c_{j},x_{i},y_{i})-g_{N}(c_{N},x_{i},y_{i})\right)
a_{j}+\sum_{j=N+1}^{j=N+N_{b}}\left(
g_{j}(c_{j},x_{i},y_{i})-g_{N}(c_{N},x_{i},y_{i})\right) \psi
_{j}(a^{1},\alpha )\text{,.}  \label{Gamma_def}
\end{equation}
Finally, substituting $a^{1}=\Gamma ^{-1}(U)$ into (\ref{fi_eq}), we arrive
at the finite dimensional continuation problem 
\begin{equation}
G(U,\alpha )=0,\;\;\text{ }U,\text{{}}G(\cdot )\in \Bbb{R}^{n\times N},\text{%
{}}\alpha \in \Bbb{R},  \label{nonlin-syst}
\end{equation}
where 
\[
G(U,\alpha )=\varphi \left( \Gamma ^{-1}(U),\text{{}}\psi \left( \Gamma
^{-1}(U),\alpha \right) ,\alpha \right) \text{,}\;\;\text{ }\Gamma :\Bbb{R}%
^{N}\rightarrow \Bbb{R}^{N}\text{, }\psi (\cdot )\in \Bbb{R}^{n\times N_{b}}.
\]

\begin{remark}
Note that in the case that the boundary condition (\ref{exact_BC}) is
linear, $\psi _{j}$ are linear, and consequently $\Gamma $ is an $N\times N$
matrix.
\end{remark}

In Section \ref{examples} we consider examples of continuation of 1D PDEs
with $\Omega ={(0,1)} $ and 2D\ PDEs with $\Omega ={(0,1)\times (0,1)} $. In
all 2D examples we have the same number $N_{s} $ of nodes in $x $ and $y $
directions. We choose a constant shape parameter $c_{j}=s/(N_{s}-1) $. Our
typical choice for $s $ is $4\leq s\leq 12 $.

We use two types of node distributions. In the case of uniform node
distribution $(x_{k}, $ $y_{l})=(kh, $ $lh) $, $k, $ $l=0,...,N_{s} $, $h=%
\frac{1}{N_{s}} $. In the case of nonuniform node distribution, the nodes
adjacent to the boundary $\partial \Omega $ are placed at the distance $%
\tilde{h}=h_{1}h $ from $\partial \Omega $, $0.1\leq h_{1}\leq 0.5 $, while
the remaining nodes are distributed uniformly. A criteria for the choice of $%
h_{1} $ was a minimum of $L_{2} $-norm of the residual in $\Omega $.

\section{Numerical experiments for 1-D and 2-D elliptic PDEs}

\label{examples}

We present several examples of continuation of solutions to systems of
nonlinear 1D and 2D elliptic PDEs. Each problem is discretized by the MQ
method described in Section \ref{MQ-method}. We then perform continuation of
the resulting system of algebraic equations (\ref{nonlin-syst}) with \textsc{%
auto97}. The principal goal of our examples is to assess the accuracy of the
detection of bifurcation points. We compare our results with some published
results and, in one case, the results of our computations with an example in 
\textsc{auto97} and \textsc{content}. We will use throughout the notation $K$
for the number of unknowns in a particular method. For our MQ method $%
K=n\times N$, where $n$ is the dimension of the system and $N$ is the number
of internal nodes. We denote by MQ(u) and MQ(nu) our MQ method with the
uniform and nonuniform node distribution, respectively.

\begin{example}
1D Gelfand-Bratu equation. This is a scalar problem 
\begin{eqnarray}
u^{^{\prime \prime }}+\lambda e^{u}{} &=&{0,\ \;\;}\text{in }\Omega ={(0,1),}
\label{1D-Bratu} \\
u{(0)} &=&u(1)=0,  \nonumber
\end{eqnarray}
that appears in combustion theory and is used as the demo example {\texttt{%
exp}} in {\textsc{auto97}} \cite{DCFKSWF97} (forth order adaptive orthogonal
spline collocation method) and demo example in {\texttt{brg}} in {\textsc{%
content}} \cite{KuzLev98} (third order adaptive finite difference method).
There is a limit (fold) point on the solution curve. We take the value of $%
\lambda $ at the limit point found from demo {\texttt{exp }}($K\geq 50$) as
exact. The following table \ref{tab:gb1d} shows comparison between numerical
results in \cite{Davidson97}, our numerical results and our experiments with 
{\textsc{content}}.
\end{example}

\begin{table}[tbp] \centering%
\caption{1D Gelfand-Bratu equation: limit point comparison} \label{tab:gb1d}
\noindent a) results for MQ correspond to uniform node distribution\newline
\begin{tabular}{|c|c|c|c|c|c|}
\hline
& \cite{DCFKSWF97}, exact & \cite{Davidson97}, $K=800$ & MQ(u), $K=5$ & 
MQ(u), $K=7$ & MQ(u), $K=9$ \\ \hline
$\lambda $ & $3.513831$ & $3.5137$ & $3.512609$ & $3.514224$ & $3.514047$ \\ 
\hline
rel. error &  & $3.7\times 10^{-5}$ & $3.5\times 10^{-4}$ & $-1.1\times
10^{-4}$ & $-6.1\times 10^{-5}$ \\ \hline
\end{tabular}

\noindent \medskip

b) results for MQ correspond to nonuniform node distribution\newline
\begin{tabular}{|c|c|c|c|c|c|}
\hline
& \cite{KuzLev98}, $K=50 $ & \cite{KuzLev98}, $K=500 $ & MQ(nu), $K=5 $ & 
MQ(nu), $K=7 $ & MQ(nu), $K=9 $ \\ \hline
$\lambda $ & $3.51145 $ & $3.51380 $ & $3.514010 $ & $3.513809 $ & $3.513828 
$ \\ \hline
rel. error & $6.8\times 10^{-4} $ & $8.\, 8\times 10^{-6} $ & $-5.1\times
10^{-5} $ & $6.\, 3\times 10^{-6} $ & $8.\, 5\times 10^{-7} $ \\ \hline
\end{tabular}
\end{table}%

\begin{example}
1D Brusselator problem. This is a reaction diffusion model for a
trimolecular chemical reaction. 
\begin{equation}
\begin{array}{l}
\frac{d_{1}}{l^{2}}{u}^{^{\prime \prime }}-(b+1)u+u^{2}v+a={0}, \\ 
u{(0)}=u(1)=a,
\end{array}
\begin{array}{l}
\frac{d_{2}}{l^{2}}{v}^{^{\prime \prime }}+bu-u^{2}v={0,\ \;\;}\text{in }%
\Omega ={(0,1),} \\ 
v{(0)}=v(1)=\frac{b}{a}.
\end{array}
\label{1D-Brusselator}
\end{equation}
A stationary bifurcation occurs \cite[Eq. (24)]{ChShMe97} at 
\[
b_{n}=1+\frac{d_{1}}{d_{2}}a^{2}+\frac{\pi ^{2}n^{2}}{l^{2}}d_{1}+\frac{l^{2}%
}{\pi ^{2}n^{2}}\frac{a^{2}}{d_{2}}>0.
\]
For $l=d_{1}=1$, $d_{2}=2$, $a=4$, $n=1,$ $2$ this gives simple
bifurcations: $b_{1}=9+\pi ^{2}+\frac{8}{\pi ^{2}}=19.680174$, $b_{2}=9+4\pi
^{2}+\frac{2}{\pi ^{2}}=48.681060$, correspondingly. For the second order
central difference method with uniform mesh of $41$ mesh points ($K=80$
unknowns), the corresponding approximate bifurcation points were found in 
\cite[Section 6.1]{ChShMe97}. The following table \ref{tab:br1d} shows
comparison between analytical, numerical results \cite[Section 6.1]{ChShMe97}
and our numerical results for values of $b_{1}$ and $b_{2}$ at simple
bifurcation points.
\end{example}

\begin{table}[tbp] \centering%

\caption{1D Brusselator equation: bifurcation points comparison} \label%
{tab:br1d} \noindent a) bifurcation point $b_{1}$\newline
\begin{tabular}{|c|c|c|c|c|c|}
\hline
& exact & \cite{ChShMe97}, $K=80$ & MQ(u), $K=10$ & MQ(u), $K=14$ & MQ(u), $%
K=18$ \\ \hline
$b_{1}$ & $19.680174$ & $19.67547$ & $19.67366$ & $19.67786$ & $19.67919$ \\ 
\hline
rel. error &  & $2.4\times 10^{-4}$ & $3.3\times 10^{-4}$ & $1.2\times
10^{-4}$ & $5.0\times 10^{-5}$ \\ \hline
\end{tabular}

\noindent \medskip

b) bifurcation point $b_{2} $\newline
\begin{tabular}{|c|c|c|c|c|c|}
\hline
& exact & \cite{ChShMe97}, $K=80 $ & MQ(u), $K=10 $ & MQ(u), $K=14 $ & 
MQ(u), $K=18 $ \\ \hline
$b_{2} $ & $48.681060 $ & $48.\, 6004 $ & $48.57476 $ & $48.63168 $ & $%
48.65605 $ \\ \hline
rel. error &  & $1.7\times 10^{-3} $ & $2.\, 2\times 10^{-3} $ & $1.0\times
10^{-3} $ & $5.\, 1\times 10^{-4} $ \\ \hline
\end{tabular}
\end{table}%

\begin{example}
Pattern formation in a 1D system with mixed boundary conditions{\textbf{\ }}%
\cite{DiMaOth1994}. 
\begin{equation}
\begin{array}{rl}
\frac{d_{1}}{\omega l^{2}}u^{^{\prime \prime }}+\beta -\kappa
u-uv^{2}=0,\;\;\;\; & \delta \frac{d_{1}}{\omega l^{2}}v^{^{\prime \prime
}}+\kappa u+uv^{2}-v=0,\;{\ \;\;}\text{in }\Omega ={(0,1)} \\ 
\theta _{1}\frac{\partial u}{\partial n}=\rho (1-\theta _{1})(\theta
_{3}u^{s}-u),\;\; & \delta \theta _{2}\frac{\partial u}{\partial n}=\delta
\rho (1-\theta _{2})(\theta _{3}v^{s}-v),\;\text{on }\partial \Omega
=\left\{ 0,1\right\} .
\end{array}
\label{Oth_reac_dif}
\end{equation}
Here $\theta _{i}\in {[0,1]},~i=1,2,3,$ are homotopy parameters. For $%
d_{1}=10^{-5}$, $\omega =10^{-2}$, $\delta =0.14$, $\beta =1.0$, $\kappa
=0.001$, $(\theta _{1},\theta _{2},\theta _{3})=(1,1,0)$ (Neumann problem).
Eq. (\ref{Oth_reac_dif}) was discretized by the second order central
difference method with equidistant mesh of $41$ mesh points ($K=80$
unknowns). The following table \cite[Table 1]{DiMaOth1994} shows a
comparison between analytic and numerical results for values of $l$ at
simple bifurcation points. 
\end{example}

\begin{table}[tbp] \centering%
\caption{1D pattern formation problem, bifurcation points} \label{tab:pf1d} 
\begin{tabular}{|c|c|c|c|c|c|c|c|c|c|}
\hline
\cite[numerical]{DiMaOth1994} & $0.047$ & $0.080$ & $0.093$ & $0.159$ & $%
0.140$ & $0.238$ & $0.186$ & $0.317$ & $0.232$ \\ \hline
\cite[analytic]{DiMaOth1994} & $0.0465$ & $0.0793$ & $0.093$ & $0.159$ & $%
0.140$ & $0.238$ & $0.186$ & $0.317$ & $0.233$ \\ \hline
MQ(nu) & $0.0465$ & $0.0793$ & $0.093$ & $0.159$ & $0.140$ & $0.238$ & $0.186
$ & $0.317$ & $0.233$ \\ \hline
\end{tabular}
\end{table}%

\medskip

Our numerical results (MQ(nu) method) with $K=18 $, coincide with the
analytic results above. In addition, we found a bifurcation point at $%
l=0.279. $

\begin{example}
2D Bratu problem 
\begin{eqnarray}
{\Delta }u+\lambda e^{u}\text{, }\Omega  &=&{(0,1)\times (0,1),}
\label{2D-Bratu} \\
u &\mid &_{\partial \Omega }=0.  \nonumber
\end{eqnarray}
This problem was studied in \cite{SchTimLos96}. It was discretized with the
second order central difference method with uniform mesh and then continued
using Implicit Block Elimination based on Recursive Projections. A limit
point was detected for some value of $\lambda $ (not reported in the paper),
and spurious limit points were detected for $K=961$, $1521,$ $2209,$ $3025$
and $\lambda $ sufficiently small. We reproduced the bifurcation diagram in 
\cite{SchTimLos96}, no spurious limit points were detected. The following
table \ref{tab:gb2d} gives the values of $\lambda $ at the limit point
computed by MQ method. 
\end{example}

\begin{table}[tbp] \centering%
\caption{2D Bratu equation, limit point} \label{tab:gb2d} \noindent 
\begin{tabular}{|c|c|c|c|c|}
\hline
& \cite{SchTimLos96}, $225\leq K\leq 3025 $ & MQ(u), $K=25 $ & MQ(u), $K=49 $
& MQ(u), $K=81 $ \\ \hline
$\lambda $ & not reported & $6.873498 $ & $6.840836 $ & $6.827400 $ \\ \hline
\end{tabular}
\end{table}%

\begin{example}
\label{2D Brus}2D Brusselator problem. 
\begin{equation}
\begin{array}{l}
\frac{d_{1}}{l^{2}}{\Delta u}-(b+1)u+u^{2}v+a={0},\;\; \\ 
u\mid _{\partial \Omega }=a,
\end{array}
\begin{array}{l}
\frac{d_{2}}{l^{2}}{\Delta v}+bu-u^{2}v={0,\ \;\;}\text{in }\Omega ={%
(0,1)\times (0,1),} \\ 
v\mid _{\partial \Omega }=\frac{b}{a}.
\end{array}
\label{2D-Brusselator}
\end{equation}
A stationary bifurcation occurs \cite[Eq. (2.26)]{ChiChe98} for 
\[
b_{m,n}=1+\frac{d_{1}}{d_{2}}a^{2}+d_{1}\pi ^{2}\left( \frac{m^{2}}{l^{2}}%
+n^{2}\right) +\frac{a^{2}}{\pi ^{2}d_{2}}\left( \frac{l^{2}}{%
m^{2}+l^{2}n^{2}}\right) >0.
\]
For $l=d_{1}=1$, $d_{2}=2$, $a=4$, $(m,n)=(1,1)$, $(m,n)=(2,2)$ this gives
simple bifurcations: $b_{1,1}=9+2\pi ^{2}+\frac{4}{\pi ^{2}}$, $%
b_{2,2}=9+8\pi ^{2}+\frac{1}{\pi ^{2}}$, correspondingly. For the second
order central difference method with equidistant mesh of $21$ mesh points,
the corresponding approximate bifurcation points are found in 
\cite[Section 5]{ChiChe98}. The following tables \ref{tab:br2d}, \ref
{tab:br2dnu} show comparisons between analytical, numerical results 
\cite[Eq. (2.26)]{ChiChe98} and our numerical results for values of $b_{1,1}$
and $b_{2,2}$ at simple bifurcation points.
\end{example}

\begin{table}[tbp] \centering%
\caption{2D Brusselator equation: bifurcation points, uniform node
distribution for MQ} \label{tab:br2d} \noindent a) bifurcation point $b_{1} $%
\newline
\begin{tabular}{|c|c|c|c|c|c|}
\hline
& exact & \cite{ChiChe98}, $K=800 $ & MQ(u), $K=50 $ & MQ(u), $K=72 $ & 
MQ(u), $K=98 $ \\ \hline
$b_{1,1} $ & $29.144494 $ & $29.104774 $ & $29.16280 $ & $29.17050 $ & $%
29.16062 $ \\ \hline
rel. error &  & $1.\, 4\times 10^{-3} $ & $-6.\, 3\times 10^{-4} $ & $-8.\,
9\times 10^{-4} $ & $-5.\, 5\times 10^{-4} $ \\ \hline
\end{tabular}

\noindent \medskip

b) bifurcation point $b_{2} $\newline
\begin{tabular}{|c|c|c|c|c|c|}
\hline
& exact & \cite{ChiChe98}, $K=800 $ & MQ(u), $K=50 $ & MQ(u), $K=72 $ & 
MQ(u), $K=98 $ \\ \hline
$b_{2,2} $ & $88.058156 $ & $87.47325 $ & $87.61578 $ & $87.86924 $ & $%
88.00143 $ \\ \hline
rel. error &  & $6.6\times 10^{-3} $ & $5.0\times 10^{-3} $ & $2.1\times
10^{-3} $ & $6.4\times 10^{-4} $ \\ \hline
\end{tabular}
\end{table}%

\medskip 
\begin{table}[tbp] \centering%
\caption{2D Brusselator equation: bifurcation points, nonuniform node
distribution for MQ} \label{tab:br2dnu} \noindent a) bifurcation point $%
b_{1} $\newline
\begin{tabular}{|c|c|c|c|c|}
\hline
& exact & MQ(nu), $K=50 $ & MQ(nu), $K=72 $ & MQ(nu), $K=98 $ \\ \hline
$b_{1,1} $ & $29.144494 $ & $29.14621 $ & $29.14726 $ & $29.14431 $ \\ \hline
rel. error &  & $-5.\, 9\times 10^{-5} $ & $-9.\, 5\times 10^{-5} $ & $6.\,
3\times 10^{-6} $ \\ \hline
\end{tabular}

\noindent \medskip

b) bifurcation point $b_{2} $\newline
\begin{tabular}{|c|c|c|c|c|}
\hline
& exact & MQ(nu), $K=50 $ & MQ(nu), $K=72 $ & MQ(nu), $K=98 $ \\ \hline
$b_{2,2} $ & $88.058156 $ & $88.15470 $ & $87.93391 $ & $88.07288 $ \\ \hline
rel. error &  & $-1.1\times 10^{-3} $ & $1.\, 4\times 10^{-3} $ & $-1.\,
7\times 10^{-4} $ \\ \hline
\end{tabular}
\end{table}%

\medskip

A Hopf bifurcation occurs \cite[Eq. (2.26)]{ChiChe98} for 
\[
b_{m,n}=1+a^{2}+(d_{1}+d_{2})\left( \frac{m^{2}}{l^{2}}+n^{2}\right) \pi ^{2}
\]
for some $m,$ $n$, and $l$ large enough. For $l=10$, $d_{1}=d_{2}=1$, $a=10$%
, $(m,n)=(1,2)$, this gives a Hopf bifurcation at $b_{1,2}=101+2\left( \frac{%
1}{100}+2^{2}\right) \pi ^{2}=180.15$\medskip , see table \ref{tab:br2dhopf}.

\begin{table}[tbp] \centering%
\caption{2D Brusselator equation, Hopf bifurcation point} \label%
{tab:br2dhopf} 
\begin{tabular}{|c|c|c|c|c|c|}
\hline
& exact & MQ(u), $K=50$ & MQ(nu), $K=50$ & MQ(u), $K=72$ & MQ(u), $K=98$ \\ 
\hline
$b_{1,2}$ & $180.\,15$ & $181.8625$ & $180.7880$ & $181.0696$ & $180.492$ \\ 
\hline
rel. error &  & $-9.5\times 10^{-3}$ & $-3.\,5\times 10^{-3}$ & $-5.1\times
10^{-3}$ & $-1.9\times 10^{-3}$ \\ \hline
\end{tabular}
\end{table}%

\section{Conclusions.\label{conclusions}}

We presented the results of our experiments with the MQ method for
continuation of solution of nonlinear 1D and 2D elliptic PDEs. We used small
number of unknowns and obtained a high accuracy for detecting bifurcation
points in our examples. Here are some sample results.

(i) For the limit point in the 1D\ Gelfand-Bratu equation, the MQ method
with $9$ unknowns gives the relative errors $6.1\times 10^{-5}$ and $%
8.5\times 10^{-7}$ for the uniform and nonuniform node distributions,
respectively. The relative error in the finite difference method with $500$
nodes is $8.8\times 10^{-6}$.

(ii) For the two bifurcation points in the 2D Brusselator problem, the MQ
method with $98$ unknowns gives the relative errors $5.5\times 10^{-4},$ $%
6.4\times 10^{-4}$ for the uniform node distribution and $6.3\times 10^{-6},$
$1.7\times 10^{-5}$ for the nonuniform node distribution. The corresponding
relative errors in the finite difference method with $800$ nodes are $%
1.4\times 10^{-3}$, $6.6\times 10^{-3}$. 

(iii)  for the first in the eigenvalue problem for the 1D Laplace operator
with $9$ unknowns gives the relative error $6\times 10^{-5}$ and $9\times
10^{-7}$ for the uniform and nonuniform node distributions, respectively.
This is equivalent in accuracy to $117$ and $950$ node solution,
respectively by the finite difference method.

The increase of the number of unknowns results in a better accuracy but also
in a larger condition number of the operator $\Gamma $ mapping solution
nodal values to the expansion coefficients. This condition number is a
limiting factor in our experiments. In the future, we plan to implement the
ideas of Kansa et al. \cite{KaHo98} to circumvent this ill conditioning
problem.

In addition we found that even a simple adaptation of the nodes adjacent to
the boundary can lead to a dramatic improvement of the accuracy in detecting
bifurcation points. Adaptive choice of the shape parameter is another way to
improve the accuracy that we plan to investigate.

Our results show that MQ method is an efficient method for continuation of
solution nonlinear PDEs.

\textbf{Acknowledgments.} Support from the NASA grant NAG8-1229 is
acknowledged by the first author.

\end{document}